\documentclass[12pt]{article}
\usepackage{amssymb}
\usepackage{amsmath}
\usepackage{amsthm}
\begin{document}
\def\Drb{\nobreak\hfill\nobreak\lower 1pt\hbox{$\square $}}
\def\black{\vrule height 6pt width 6pt depth 1pt} 
\newcommand{\Aut}{\mbox{Aut }}
\newcommand{\shift}{\mbox{shift}}
\newcommand{\Ad}{\mbox{Ad }}
\newcommand{\Fun}{\mbox{Fun }}
\newcommand{\ind}{\mbox{ind }}
\newcommand{\diag}{\mbox{diag }}
\newcommand{\Ind}{\mbox{Ind }}
\newcommand{\spann}{\mbox{span }}
\newcommand{\Ob}{\mbox{Ob }}
\newcommand{\Hom}{\mbox{Hom }}
\newcommand{\pfeil}{\longrightarrow}
\newcommand{\Pfeil}{\longmapsto}
\newcommand{\proj}{\mbox{proj }}
\renewcommand{\blacksquare}{\nobreak \hfill \nobreak \black}
\newcommand{\blacksquared}{\nobreak \hfill \nobreak \black}
\newcommand{\tm}{\otimes _M}
\newcommand{\tn}{\otimes _N} 
\newcommand{\tp}{\otimes _P}
\newcommand{\tq}{\otimes _Q}
\newcommand{\ta}{\otimes _A}
\renewcommand{\th}{\otimes _H}
\newcommand{\kreuz}{\rtimes}
\newcommand{\ckreuz}{{\rtimes} _c}
\newcommand{\ockreuz}{{\rtimes} _{\bar{c}}}
\newcommand{\rest}{\! \mid \!}
\newcommand{\C}[1]{{\cal #1}}
\newcommand{\B}[1]{{\bf #1}}
\newcommand{\equi}{\Longleftrightarrow}
\newcommand{\beq}{\begin{equation}}
\newcommand{\eeq}{\end{equation}}
\newcommand{\Kl}[3]{\langle #1 \otimes #2, #3\rangle}
\newcommand{\ska}[2]{\langle #1,\, #2\rangle}
\newcommand{\skab}[2]{\Bigl\langle #1,\, #2\Bigr\rangle}
\newcommand{\skal}[2]{\langle #1,\, #2\rangle _N^l}
\newcommand{\skar}[2]{\langle #1,\, #2\rangle _N^r}
\newcommand{\ov}{\overline}
\newcommand{\newl}{\newline \noindent}
\newcommand{\tr}{\mbox{tr}}
\newcommand{\Tr}{\mbox{Tr}}
\newcommand{\zws}{{\rm II}_1}
\newcommand{\abs}{\vspace{4mm}}
\newcommand{\orho}{\bar{\rho }}
\newcommand{\osigma}{\bar{\sigma }}
\newcommand{\ophi}{\bar{\phi }}
\newcommand{\opsi}{\bar{\psi }}
\newcommand{\otau}{\bar{\tau }}
\newcommand{\ovs}[1]{\ov{\sigma (#1)}}
\newcommand{\oR}{\bar{R}}
\newcommand{\oV}{\bar{V}}
\newcommand{\oW}{\bar{W}}
\newcommand{\oc}{\bar{c}}
\newcommand{\comp}{\mathbb{C}}
\newcommand{\nat}{\mathbb{N}}
\newcommand{\reel}{\mathbb{R}}
\newcommand{\ganz}{\mathbb{Z}}
\newcommand{\rat}{\mathbb{Q}}

 \title{$\zws $-Subfactors associated with the $C^*$-Tensor Category
of a Finite Group}
\author{Reinhard Schaf\-litzel\\
        Mathematisches Institut\\
        Technische Universit\"at M\"unchen\\
        Arcisstr. 21\\ 
        80290 M\"unchen\\
        GERMANY\\
        E-mail: schafl@mathematik.tu-muenchen.de}
\date {December 2, 1996}

\maketitle

\theoremstyle{plain}
\newtheorem{lem}[subsection]{Lemma}
\newtheorem{cor}[subsection]{Corollary}
\newtheorem{prop}[subsection]{Proposition}
\newtheorem{thm}[subsection]{Theorem}
\theoremstyle{definition}
\newtheorem{defi}[subsection]{Definition}
\newtheorem{quest}[subsection]{Question}
\newtheorem{ass}[subsection]{Assumption}
\setcounter{tocdepth}{1}
\hyphenation{bundles Yamagami}
%
%
%
 
\long\def\ig#1{\relax}
\ig{Thanks to Roberto Minio for this def'n.  Compare the def'n of
\comment in AMSTeX.}
 
\setlength{\unitlength}{.01em}%
 
\newcount \coefa
\newcount \coefb
\newcount \coefc
\newdimen\tempdimen
\newdimen\xlen
\newdimen\ylen
 
\newcount\tempcounta
\newcount\tempcountb
\newcount\tempcountc
\newcount\tempcountd
\newcount\tempcounte
\newcount\tempcountf
\newcount\xext
\newcount\yext
\newcount\xoff
\newcount\yoff
\newcount\gap%
\newcount\arrowtypea
\newcount\arrowtypeb
\newcount\arrowtypec
\newcount\arrowtyped
\newcount\arrowtypee
\newcount\height
\newcount\width
\newcount\xpos
\newcount\ypos
\newcount\run
\newcount\rise
\newcount\arrowlength
\newcount\halflength
\newcount\arrowtype
\newsavebox{\tempboxa}%
\newsavebox{\tempboxb}%
\newsavebox{\tempboxc}%

\catcode`@=11 
\def\settoheight#1#2{\setbox\@tempboxa\hbox{#2}#1\ht\@tempboxa\relax}%
\def\settodepth#1#2{\setbox\@tempboxa\hbox{#2}#1\dp\@tempboxa\relax}%
\let\ifnextchar=\@ifnextchar
\catcode`@=12 
 
\def\putbox(#1,#2)#3{\put(#1,#2){\makebox(0,0){#3}}}

\def\setsqparms[#1`#2`#3`#4;#5`#6]{%
\settripairparms[#1`#2`#3`#4`1;#6]%
\width #5
}
 
\def\settriparms[#1`#2`#3;#4]{\settripairparms[#1`#2`#3`1`1;#4]}%

\def\settripairparms[#1`#2`#3`#4`#5;#6]{%
\arrowtypea #1
\arrowtypeb #2
\arrowtypec #3
\arrowtyped #4
\arrowtypee #5
\height #6
\width #6
}
 
\def\resetparms{\settripairparms[1`1`1`1`1;500]\width 500}
 
\def\mvector(#1,#2)#3{
\put(0,0){\vector(#1,#2){#3}}%
\put(0,0){\vector(#1,#2){30}}%
}
\def\evector(#1,#2)#3{{
\arrowlength #3
\put(0,0){\vector(#1,#2){\arrowlength}}%
\advance \arrowlength by-30
\put(0,0){\vector(#1,#2){\arrowlength}}%
}}

\def\horsize#1#2{%
\settowidth{\tempdimen}{$#2$}%
#1=\tempdimen
\divide #1 by\unitlength
}
 
\def\vertsize#1#2{%
\settoheight{\tempdimen}{$#2$}%
#1=\tempdimen
\settodepth{\tempdimen}{$#2$}%
\advance #1 by\tempdimen
\divide #1 by\unitlength
}

\def\vertadjust[#1`#2`#3]{%
\vertsize{\tempcounta}{#1}%
\vertsize{\tempcountb}{#2}%
\ifnum \tempcounta<\tempcountb \tempcounta=\tempcountb \fi
\divide\tempcounta by2
\vertsize{\tempcountb}{#3}%
\ifnum \tempcountb>0 \advance \tempcountb by20 \fi
\ifnum \tempcounta<\tempcountb \tempcounta=\tempcountb \fi
}
 
\def\horadjust[#1`#2`#3]{%
\horsize{\tempcounta}{#1}%
\horsize{\tempcountb}{#2}%
\ifnum \tempcounta<\tempcountb \tempcounta=\tempcountb \fi
\divide\tempcounta by2
\horsize{\tempcountb}{#3}%
\ifnum \tempcountb>0 \advance \tempcountb by60 \fi
\ifnum \tempcounta<\tempcountb \tempcounta=\tempcountb \fi
}
 
\ig{ In this procedure, #1 is the paramater that sticks out all the way,
#2 sticks out the least and #3 is a label sticking out half way.  #4 is
the amount of the offset.}
 
\def\sladjust[#1`#2`#3]#4{%
\tempcountc=#4
\horsize{\tempcounta}{#1}%
\divide \tempcounta by2
\horsize{\tempcountb}{#2}%
\divide \tempcountb by2
\advance \tempcountb by-\tempcountc
\ifnum \tempcounta<\tempcountb \tempcounta=\tempcountb\fi
\divide \tempcountc by2
\horsize{\tempcountb}{#3}%
\advance \tempcountb by-\tempcountc
\ifnum \tempcountb>0 \advance \tempcountb by80\fi
\ifnum \tempcounta<\tempcountb \tempcounta=\tempcountb\fi
\advance\tempcounta by20
}
 
\def\putvector(#1,#2)(#3,#4)#5#6{{%
\xpos=#1
\ypos=#2
\run=#3
\rise=#4
\arrowlength=#5
\arrowtype=#6
\ifnum \arrowtype<0
    \ifnum \run=0
        \advance \ypos by-\arrowlength
    \else
        \tempcounta \arrowlength
        \multiply \tempcounta by\rise
        \divide \tempcounta by\run
        \ifnum\run>0
            \advance \xpos by\arrowlength
            \advance \ypos by\tempcounta
        \else
            \advance \xpos by-\arrowlength
            \advance \ypos by-\tempcounta
        \fi
    \fi
    \multiply \arrowtype by-1
    \multiply \rise by-1
    \multiply \run by-1
\fi
\ifnum \arrowtype=1
    \put(\xpos,\ypos){\vector(\run,\rise){\arrowlength}}%
\else\ifnum \arrowtype=2
    \put(\xpos,\ypos){\mvector(\run,\rise)\arrowlength}%
\else\ifnum\arrowtype=3
    \put(\xpos,\ypos){\evector(\run,\rise){\arrowlength}}%
\fi\fi\fi
}}
 
\def\bfig{\begin{picture}(\xext,\yext)(\xoff,\yoff)}
\def\efig{\end{picture}}

\def\putsplitvector(#1,#2)#3#4{
\xpos #1
\ypos #2
\arrowtype #4
\halflength #3
\arrowlength #3
\gap 140
\advance \halflength by-\gap
\divide \halflength by2
\ifnum \arrowtype=1
    \put(\xpos,\ypos){\line(0,-1){\halflength}}%
    \advance\ypos by-\halflength
    \advance\ypos by-\gap
    \put(\xpos,\ypos){\vector(0,-1){\halflength}}%
\else\ifnum \arrowtype=2
    \put(\xpos,\ypos){\line(0,-1)\halflength}%
    \put(\xpos,\ypos){\vector(0,-1)3}%
    \advance\ypos by-\halflength
    \advance\ypos by-\gap
    \put(\xpos,\ypos){\vector(0,-1){\halflength}}%
\else\ifnum\arrowtype=3
    \put(\xpos,\ypos){\line(0,-1)\halflength}%
    \advance\ypos by-\halflength
    \advance\ypos by-\gap
    \put(\xpos,\ypos){\evector(0,-1){\halflength}}%
\else\ifnum \arrowtype=-1
    \advance \ypos by-\arrowlength
    \put(\xpos,\ypos){\line(0,1){\halflength}}%
    \advance\ypos by\halflength
    \advance\ypos by\gap
    \put(\xpos,\ypos){\vector(0,1){\halflength}}%
\else\ifnum \arrowtype=-2
    \advance \ypos by-\arrowlength
    \put(\xpos,\ypos){\line(0,1)\halflength}%
    \put(\xpos,\ypos){\vector(0,1)3}%
    \advance\ypos by\halflength
    \advance\ypos by\gap
    \put(\xpos,\ypos){\vector(0,1){\halflength}}%
\else\ifnum\arrowtype=-3
    \advance \ypos by-\arrowlength
    \put(\xpos,\ypos){\line(0,1)\halflength}%
    \advance\ypos by\halflength
    \advance\ypos by\gap
    \put(\xpos,\ypos){\evector(0,1){\halflength}}%
\fi\fi\fi\fi\fi\fi
}
 
\def\setpos(#1,#2){\xpos=#1 \ypos#2}
 
\def\putmorphism(#1)(#2,#3)[#4`#5`#6]#7#8#9{{%
\run #2
\rise #3
\ifnum\rise=0
  \puthmorphism(#1)[#4`#5`#6]{#7}{#8}{#9}%
\else\ifnum\run=0
  \putvmorphism(#1)[#4`#5`#6]{#7}{#8}{#9}%
\else
\setpos(#1)%
\arrowlength #7
\arrowtype #8
\ifnum\run=0
\else\ifnum\rise=0
\else
\ifnum\run>0
    \coefa=1
\else
   \coefa=-1
\fi
\ifnum\arrowtype>0
   \coefb=0
   \coefc=-1
\else
   \coefb=\coefa
   \coefc=1
   \arrowtype=-\arrowtype
\fi
\width=2
\multiply \width by\run
\divide \width by\rise
\ifnum \width<0  \width=-\width\fi
\advance\width by60
\if l#9 \width=-\width\fi
\putbox(\xpos,\ypos){$#4$}
{\multiply \coefa by\arrowlength
\advance\xpos by\coefa
\multiply \coefa by\rise
\divide \coefa by\run
\advance \ypos by\coefa
\putbox(\xpos,\ypos){$#5$} }%
{\multiply \coefa by\arrowlength
\divide \coefa by2
\advance \xpos by\coefa
\advance \xpos by\width
\multiply \coefa by\rise
\divide \coefa by\run
\advance \ypos by\coefa
\if l#9%
   \put(\xpos,\ypos){\makebox(0,0)[r]{$#6$}}%
\else\if r#9%
   \put(\xpos,\ypos){\makebox(0,0)[l]{$#6$}}%
\fi\fi }%
{\multiply \rise by-\coefc
\multiply \run by-\coefc
\multiply \coefb by\arrowlength
\advance \xpos by\coefb
\multiply \coefb by\rise
\divide \coefb by\run
\advance \ypos by\coefb
\multiply \coefc by70
\advance \ypos by\coefc
\multiply \coefc by\run
\divide \coefc by\rise
\advance \xpos by\coefc
\multiply \coefa by140
\multiply \coefa by\run
\divide \coefa by\rise
\advance \arrowlength by\coefa
\ifnum \arrowtype=1
   \put(\xpos,\ypos){\vector(\run,\rise){\arrowlength}}%
\else\ifnum\arrowtype=2
   \put(\xpos,\ypos){\mvector(\run,\rise){\arrowlength}}%
\else\ifnum\arrowtype=3
   \put(\xpos,\ypos){\evector(\run,\rise){\arrowlength}}%
\fi\fi\fi}%
\fi\fi
\fi\fi}}

\def\puthmorphism(#1,#2)[#3`#4`#5]#6#7#8{{%
\xpos #1
\ypos #2
\width #6
\arrowlength #6
\putbox(\xpos,\ypos){$#3$\vphantom{$#4$}}%
{\advance \xpos by\arrowlength
\putbox(\xpos,\ypos){\vphantom{$#3$}$#4$}}%
\horsize{\tempcounta}{#3}%
\horsize{\tempcountb}{#4}%
\divide \tempcounta by2
\divide \tempcountb by2
\advance \tempcounta by30
\advance \tempcountb by30
\advance \xpos by\tempcounta
\advance \arrowlength by-\tempcounta
\advance \arrowlength by-\tempcountb
\putvector(\xpos,\ypos)(1,0){\arrowlength}{#7}%
\divide \arrowlength by2
\advance \xpos by\arrowlength
\vertsize{\tempcounta}{#5}%
\divide\tempcounta by2
\advance \tempcounta by20
\if a#8 %
   \advance \ypos by\tempcounta
   \put(\xpos,\ypos){\makebox(0,0){$#5$}}%
\else
   \advance \ypos by-\tempcounta
   \put(\xpos,\ypos){\makebox(0,0){$#5$}}%
\fi
}}
 
\def\putvmorphism(#1,#2)[#3`#4`#5]#6#7#8{{%
\xpos #1
\ypos #2
\arrowlength #6
\arrowtype #7
\settowidth{\xlen}{$#5$}%
\putbox(\xpos,\ypos){$#3$}%
{\advance \ypos by-\arrowlength
\putbox(\xpos,\ypos){$#4$}}%
{\advance\arrowlength by-140
\advance \ypos by -70
\ifdim\xlen>0pt
   \if m#8%
      \putsplitvector(\xpos,\ypos){\arrowlength}{\arrowtype}%
   \else
      \putvector(\xpos,\ypos)(0,-1){\arrowlength}{\arrowtype}%
   \fi
\else
   \putvector(\xpos,\ypos)(0,-1){\arrowlength}{\arrowtype}%
\fi}%
\ifdim\xlen>0pt
   \divide \arrowlength by2
   \advance\ypos by-\arrowlength
   \if l#8%
      \advance \xpos by-40
      \put(\xpos,\ypos){\makebox(0,0)[r]{$#5$}}%
   \else\if r#8%
      \advance \xpos by40
      \put(\xpos,\ypos){\makebox(0,0)[l]{$#5$}}%
   \else
      \putbox(\xpos,\ypos){$#5$}%
   \fi\fi
\fi
}}
 
\def\topadjust[#1`#2`#3]{%
\yoff=10
\vertadjust[#1`#2`{#3}]%
\advance \yext by\tempcounta
\advance \yext by 10
}
 
\def\botadjust[#1`#2`#3]{%
\vertadjust[#1`#2`{#3}]%
\advance \yext by\tempcounta
\advance \yoff by-\tempcounta
}
 
\def\leftadjust[#1`#2`#3]{%
\xoff=0
\horadjust[#1`#2`{#3}]%
\advance \xext by\tempcounta
\advance \xoff by-\tempcounta
}
 
\def\rightadjust[#1`#2`#3]{%
\horadjust[#1`#2`{#3}]%
\advance \xext by\tempcounta
}
 
\def\rightsladjust[#1`#2`#3]{%
\sladjust[#1`#2`{#3}]{\width}%
\advance \xext by\tempcounta
}
 
\def\leftsladjust[#1`#2`#3]{%
\xoff=0
\sladjust[#1`#2`{#3}]{\width}%
\advance \xext by\tempcounta
\advance \xoff by-\tempcounta
}
 
\def\adjust[#1`#2;#3`#4;#5`#6;#7`#8]{%
\topadjust[#1``{#2}]
\leftadjust[#3``{#4}]
\rightadjust[#5``{#6}]
\botadjust[#7``{#8}]}

\def\putsquare(#1)[#2`#3`#4`#5;#6`#7`#8`#9]{%
\setpos(#1)
\puthmorphism(\xpos,\ypos)[#4`#5`{#9}]{\width}{\arrowtyped}b%
\advance\ypos by \height
\puthmorphism(\xpos,\ypos)[#2`#3`{#6}]{\width}{\arrowtypea}a%
\putvmorphism(\xpos,\ypos)[``{#7}]{\height}{\arrowtypeb}l%
\advance\xpos by \width
\putvmorphism(\xpos,\ypos)[``{#8}]{\height}{\arrowtypec}r%
}
 
\def\square[#1`#2`#3`#4;#5`#6`#7`#8]{{
\xext=\width                              
\yext=\height                             
\topadjust[#1`#2`{#5}]
\botadjust[#3`#4`{#8}]
\leftadjust[#1`#3`{#6}]
\rightadjust[#2`#4`{#7}]
\begin{picture}(\xext,\yext)(\xoff,\yoff)
\putsquare(0,0)[#1`#2`#3`#4;#5`#6`#7`{#8}]
\end{picture}
}}

\def\putptriangle(#1,#2)[#3`#4`#5;#6`#7`#8]{%
\xpos=#1 \ypos=#2
\advance\ypos by \height
\puthmorphism(\xpos,\ypos)[#3`#4`{#6}]{\height}{\arrowtypea}a%
\putvmorphism(\xpos,\ypos)[`#5`{#7}]{\height}{\arrowtypeb}l%
\advance\xpos by\height
\putmorphism(\xpos,\ypos)(-1,-1)[``{#8}]{\height}{\arrowtypec}r%
}
 
\def\ptriangle[#1`#2`#3;#4`#5`#6]{{
\width=\height                         
\xext=\width                           
\yext=\width                           
\topadjust[#1`#2`{#4}]
\botadjust[#3``]
\leftadjust[#1`#3`{#5}]
\rightsladjust[#2`#3`{#6}]
\begin{picture}(\xext,\yext)(\xoff,\yoff)
\putptriangle(0,0)[#1`#2`#3;#4`#5`{#6}]%
\end{picture}%
}}

\def\putqtriangle(#1,#2)[#3`#4`#5;#6`#7`#8]{%
\xpos=#1 \ypos=#2
\advance\ypos by\height
\puthmorphism(\xpos,\ypos)[#3`#4`{#6}]{\height}{\arrowtypea}a%
\putmorphism(\xpos,\ypos)(1,-1)[``{#7}]{\height}{\arrowtypeb}l%
\advance\xpos by\height
\putvmorphism(\xpos,\ypos)[`#5`{#8}]{\height}{\arrowtypec}r%
}
 
\def\qtriangle[#1`#2`#3;#4`#5`#6]{{
\width=\height                         
\xext=\width                           
\yext=\height                          
\topadjust[#1`#2`{#4}]
\botadjust[#3``]
\leftsladjust[#1`#3`{#5}]
\rightadjust[#2`#3`{#6}]
\begin{picture}(\xext,\yext)(\xoff,\yoff)
\putqtriangle(0,0)[#1`#2`#3;#4`#5`{#6}]%
\end{picture}%
}}

\def\putdtriangle(#1,#2)[#3`#4`#5;#6`#7`#8]{%
\xpos=#1 \ypos=#2
\puthmorphism(\xpos,\ypos)[#4`#5`{#8}]{\height}{\arrowtypec}b%
\advance\xpos by \height \advance\ypos by\height
\putmorphism(\xpos,\ypos)(-1,-1)[``{#6}]{\height}{\arrowtypea}l%
\putvmorphism(\xpos,\ypos)[#3``{#7}]{\height}{\arrowtypeb}r%
}
 
\def\dtriangle[#1`#2`#3;#4`#5`#6]{{
\width=\height                         
\xext=\width                           
\yext=\height                          
\topadjust[#1``]
\botadjust[#2`#3`{#6}]
\leftsladjust[#2`#1`{#4}]
\rightadjust[#1`#3`{#5}]
\begin{picture}(\xext,\yext)(\xoff,\yoff)
\putdtriangle(0,0)[#1`#2`#3;#4`#5`{#6}]%
\end{picture}%
}}

\def\putbtriangle(#1,#2)[#3`#4`#5;#6`#7`#8]{%
\xpos=#1 \ypos=#2
\puthmorphism(\xpos,\ypos)[#4`#5`{#8}]{\height}{\arrowtypec}b%
\advance\ypos by\height
\putmorphism(\xpos,\ypos)(1,-1)[``{#7}]{\height}{\arrowtypeb}r%
\putvmorphism(\xpos,\ypos)[#3``{#6}]{\height}{\arrowtypea}l%
}
 
\def\btriangle[#1`#2`#3;#4`#5`#6]{{
\width=\height                         
\xext=\width                           
\yext=\height                          
\topadjust[#1``]
\botadjust[#2`#3`{#6}]
\leftadjust[#1`#2`{#4}]
\rightsladjust[#3`#1`{#5}]
\begin{picture}(\xext,\yext)(\xoff,\yoff)
\putbtriangle(0,0)[#1`#2`#3;#4`#5`{#6}]%
\end{picture}%
}}

\def\putAtriangle(#1,#2)[#3`#4`#5;#6`#7`#8]{%
\xpos=#1 \ypos=#2
{\multiply \height by2
\puthmorphism(\xpos,\ypos)[#4`#5`{#8}]{\height}{\arrowtypec}b}%
\advance\xpos by\height \advance\ypos by\height
\putmorphism(\xpos,\ypos)(-1,-1)[#3``{#6}]{\height}{\arrowtypea}l%
\putmorphism(\xpos,\ypos)(1,-1)[``{#7}]{\height}{\arrowtypeb}r%
}
 
\def\Atriangle[#1`#2`#3;#4`#5`#6]{{
\width=\height                         
\xext=\width                           
\yext=\height                          
\topadjust[#1``]
\botadjust[#2`#3`{#6}]
\multiply \xext by2 
\leftsladjust[#2`#1`{#4}]
\rightsladjust[#3`#1`{#5}]
\begin{picture}(\xext,\yext)(\xoff,\yoff)%
\putAtriangle(0,0)[#1`#2`#3;#4`#5`{#6}]%
\end{picture}%
}}

\def\putAtrianglepair(#1,#2)[#3]{\xpos=#1 \ypos=#2%
\putAtrianglepairx[#3]}
\def\putAtrianglepairx[#1`#2`#3`#4;#5`#6`#7`#8`#9]{%
\puthmorphism(\xpos,\ypos)[#2`#3`{#8}]{\height}{\arrowtyped}b%
\advance\xpos by\height
\puthmorphism(\xpos,\ypos)[\phantom{#3}`#4`{#9}]{\height}{\arrowtypee}b%
\advance\ypos by\height
\putmorphism(\xpos,\ypos)(-1,-1)[#1``{#5}]{\height}{\arrowtypea}l%
\putvmorphism(\xpos,\ypos)[``{#6}]{\height}{\arrowtypeb}m%
\putmorphism(\xpos,\ypos)(1,-1)[``{#7}]{\height}{\arrowtypec}r%
}
 
\def\Atrianglepair[#1`#2`#3`#4;#5`#6`#7`#8`#9]{{%
\width=\height
\xext=\width
\yext=\height
\topadjust[#1``]%
\vertadjust[#2`#3`{#8}]
\tempcountd=\tempcounta                       
\vertadjust[#3`#4`{#9}]
\ifnum\tempcounta<\tempcountd                 
\tempcounta=\tempcountd\fi                    
\advance \yext by\tempcounta                  
\advance \yoff by-\tempcounta                 
\multiply \xext by2 
\leftsladjust[#2`#1`{#5}]
\rightsladjust[#4`#1`{#7}]%
\begin{picture}(\xext,\yext)(\xoff,\yoff)%
\putAtrianglepair(0,0)[#1`#2`#3`#4;#5`#6`#7`#8`{#9}]%
\end{picture}%
}}

\def\putVtriangle(#1,#2)[#3`#4`#5;#6`#7`#8]{%
\xpos=#1 \ypos=#2
\advance\ypos by\height
{\multiply\height by2
\puthmorphism(\xpos,\ypos)[#3`#4`{#6}]{\height}{\arrowtypea}a}%
\putmorphism(\xpos,\ypos)(1,-1)[`#5`{#7}]{\height}{\arrowtypeb}l%
\advance\xpos by\height
\advance\xpos by\height
\putmorphism(\xpos,\ypos)(-1,-1)[``{#8}]{\height}{\arrowtypec}r%
}
 
\def\Vtriangle[#1`#2`#3;#4`#5`#6]{{
\width=\height                         
\xext=\width                           
\yext=\height                          
\topadjust[#1`#2`{#4}]
\botadjust[#3``]
\multiply \xext by2 
\leftsladjust[#1`#3`{#5}]
\rightsladjust[#2`#3`{#6}]
\begin{picture}(\xext,\yext)(\xoff,\yoff)%
\putVtriangle(0,0)[#1`#2`#3;#4`#5`{#6}]%
\end{picture}%
}}

\def\putVtrianglepair(#1,#2)[#3]{\xpos=#1 \ypos=#2%
\putVtrianglepairx[#3]}
\def\putVtrianglepairx[#1`#2`#3`#4;#5`#6`#7`#8`#9]{%
\advance\ypos by\height
\putmorphism(\xpos,\ypos)(1,-1)[`#4`{#7}]{\height}{\arrowtypec}l%
\puthmorphism(\xpos,\ypos)[#1`#2`{#5}]{\height}{\arrowtypea}a%
\advance\xpos by\height
\puthmorphism(\xpos,\ypos)[\phantom{#2}`#3`{#6}]{\height}{\arrowtypeb}a%
\putvmorphism(\xpos,\ypos)[``{#8}]{\height}{\arrowtyped}m%
\advance\xpos by\height
\putmorphism(\xpos,\ypos)(-1,-1)[``{#9}]{\height}{\arrowtypee}r%
}
 
\def\Vtrianglepair[#1`#2`#3`#4;#5`#6`#7`#8`#9]{{%
\xoff=0
\yoff=2 
\xext=\height                  
\width=\height                 
\yext=\height                  
\vertadjust[#1`#2`{#5}]
\tempcountd=\tempcounta        
\vertadjust[#2`#3`{#6}]
\ifnum\tempcounta<\tempcountd
\tempcounta=\tempcountd\fi
\advance \yext by\tempcounta
\botadjust[#4``]%
\multiply \xext by2
\leftsladjust[#1`#4`{#7}]%
\rightsladjust[#3`#4`{#9}]%
\begin{picture}(\xext,\yext)(\xoff,\yoff)%
\putVtrianglepair(0,0)[#1`#2`#3`#4;#5`#6`#7`#8`{#9}]%
\end{picture}%
}}

\def\putCtriangle(#1,#2)[#3`#4`#5;#6`#7`#8]{%
\xpos=#1 \ypos=#2
\advance\ypos by\height
\putmorphism(\xpos,\ypos)(1,-1)[``{#8}]{\height}{\arrowtypec}l%
\advance\xpos by\height
\advance\ypos by\height
\putmorphism(\xpos,\ypos)(-1,-1)[#3`#4`{#6}]{\height}{\arrowtypea}l%
{\multiply\height by 2
\putvmorphism(\xpos,\ypos)[`#5`{#7}]{\height}{\arrowtypeb}r}%
}
 
\def\Ctriangle[#1`#2`#3;#4`#5`#6]{{
\width=\height                          
\xext=\width                            
\yext=\height                           
\multiply \yext by2 
\topadjust[#1``]
\botadjust[#3``]
\sladjust[#2`#1`{#4}]{\width}
\tempcountd=\tempcounta                 
\sladjust[#2`#3`{#6}]{\width}
\ifnum \tempcounta<\tempcountd          
\tempcounta=\tempcountd\fi              
\advance \xext by\tempcounta            
\advance \xoff by-\tempcounta           
\rightadjust[#1`#3`{#5}]
\begin{picture}(\xext,\yext)(\xoff,\yoff)%
\putCtriangle(0,0)[#1`#2`#3;#4`#5`{#6}]%
\end{picture}%
}}

\def\putDtriangle(#1,#2)[#3`#4`#5;#6`#7`#8]{%
\xpos=#1 \ypos=#2
\advance\xpos by\height \advance\ypos by\height
\putmorphism(\xpos,\ypos)(-1,-1)[``{#8}]{\height}{\arrowtypec}r%
\advance\xpos by-\height \advance\ypos by\height
\putmorphism(\xpos,\ypos)(1,-1)[`#4`{#7}]{\height}{\arrowtypeb}r%
{\multiply\height by 2
\putvmorphism(\xpos,\ypos)[#3`#5`{#6}]{\height}{\arrowtypea}l}%
}
 
\def\Dtriangle[#1`#2`#3;#4`#5`#6]{{
\width=\height                         
\xext=\height                          
\yext=\height                          
\multiply \yext by2 
\topadjust[#1``]
\botadjust[#3``]
\leftadjust[#1`#3`{#4}]
\sladjust[#2`#1`{#4}]{\height}
\tempcountd=\tempcountd                
\sladjust[#2`#3`{#6}]{\height}
\ifnum \tempcounta<\tempcountd         
\tempcounta=\tempcountd\fi             
\advance \xext by\tempcounta           
\begin{picture}(\xext,\yext)(\xoff,\yoff)
\putDtriangle(0,0)[#1`#2`#3;#4`#5`{#6}]%
\end{picture}%
}}

\def\setrecparms[#1`#2]{\width=#1 \height=#2}%
\def\recurse[#1`#2`#3`#4;#5`#6`#7`#8`#9]{{%
\settowidth{\tempdimen}{#1}
\ifdim\tempdimen=0pt
  \savebox{\tempboxa}{\hbox{#2}}%
  \savebox{\tempboxb}{\hbox{#4}}%
  \savebox{\tempboxc}{\hbox{#7}}%
\else
  \savebox{\tempboxa}{\hbox{$\hbox{#1}\times\hbox{#2}$}}%
  \savebox{\tempboxb}{\hbox{$\hbox{#1}\times\hbox{#4}$}}%
  \savebox{\tempboxc}{\hbox{$\hbox{#1}\times\hbox{#7}$}}%
\fi
\tempcounte=\height
\divide\tempcounte by 2
\tempcountf=\tempcounte
\advance\tempcountf by \width
\xext=\tempcountf \yext=\height
\topadjust[#2`\usebox{\tempboxa}`{#5}]%
\botadjust[#4`\usebox{\tempboxb}`{#9}]%
\sladjust[#3`#2`{#6}]{\tempcounte}%
\tempcountd=\tempcounta
\sladjust[#3`#4`{#8}]{\tempcounte}%
\ifnum \tempcounta<\tempcountd
\tempcounta=\tempcountd\fi
\advance \xext by\tempcounta
\advance \xoff by-\tempcounta
\rightadjust[\usebox{\tempboxa}`\usebox{\tempboxb}`\usebox{\tempboxc}]%
\bfig
{\settriparms[-1`1`1;\tempcounte]%
\putCtriangle(0,0)[`#3`;#6`#7`{#8}]}%
\arrowtypea=-1 \arrowtypeb=0 \arrowtypec=1 \arrowtyped=-1
\putsquare(\tempcounte,0)[#2`\usebox{\tempboxa}`#4`\usebox{\tempboxb};%
#5``\usebox{\tempboxc}`#9]%
\efig
}}

\begin{abstract} 
We determine the subfactors $N\subset R$ of the hyperfinite $\zws $-factor
$R$ with finite index for which the $C^*$-tensor category of the 
associated $(N,N)$-bimodules is equivalent to the $C^*$-tensor category
$\C{U}_G$ of all unitary finite dimensional representations of a given
finite group $G$. It turns out that every subfactor of that kind is
isomorphic to a subfactor $R^G\subset (R\otimes \B{L}(\comp ^r))^H$,
where $R^G$ is the fixed point algebra under an outer action $\alpha $ of
$G$, $H$ is a subgroup of $G$, $\psi :H\pfeil \B{U}(\comp ^R)$ is a 
unitary finite
dimensional projective representation of $H$ satisfying a certain additional
condition and $(R\otimes \B{L}(\comp ^r))^H$ is the fixed point algebra under 
the action $\alpha \rest H \otimes \Ad \psi $ of $H$
on $R\otimes \B{L}(\comp ^r)$.

\end{abstract}

\section*{Introduction} 
Recently the concept of $C^*$-tensor categories appeared in Jones' theory
of subfactors. Mainly we are interested in $C^*$-tensor categories as
introduced by J.E. Roberts and R. Longo in \cite{LoRo} (see also
\cite{Ro}, \cite{FK} and \cite{Sch4}), in which every object has a conjugate.
$C^*$-tensor categories of that kind are called compact in this article.
The most natural example
of a compact $C^*$-tensor category is the category $\C{U}_G$ of the unitary
finite dimensional representations of a compact group $G$. Compact
$C^*$-tensor categories may be regarded as a concept to deal with more general
symmetries than those described by groups. In particular the unitary
finite dimensional corepresentations of a compact matrix pseudogroup
(in the sense of S.L. Woronowicz \cite{Wor})  form a compact 
$C^*$-tensor category, but the concept of a compact $C^*$-tensor category
is more comprehensive than that of a compact quantum group.

The $(N,N)$-bimodules belonging to a $\zws $-subfactor $N\subset M$ with 
finite Jones index form a compact $C^*$-tensor category $\C{B}_{N\subset M}$.
The number of
the equivalence classes of the irreducible objects of $\C{B}_{N\subset M}$ is
finite if and only if the inclusion $N\subset M$ has finite depth. 
As S. Popa has shown,
the subfactors of the hyperfinite $\zws $-factor with finite depth are 
completely classified by their standard invariant, which is given by a single
commuting square of finite dimensional $*$-algebras
(see \cite{P1}). A. Ocneanu found the 
so-called flatness condition which allows to decide which of the commuting
squares are standard invariants of finite depth subfactors
(see \cite{Oc2}).
Unfortunately, checking the flatness condition requires very complicated
computations in most cases.

But using $C^*$-tensor categories we get another method for the classification
of finite depth subfactors of the hyperfinite factor. It requires the solution
of two problems:
\begin{itemize}
\item[1)] Classify all finite $C^*$-tensor categories up to equivalence
(i. e. all compact $C^*$-tensor categories for which the number
of equivalence classes of irreducible objects is finite).
\item[2)] Classify all subfactors $N\subset R$ of the hyperfinite 
$\zws $-factor $R$ with finite index, for which the associated 
$C^*$-tensor category of the $(N,N)$-bimodules is equivalent to a given
$C^*$-tensor category $\C{C}$.
\end{itemize}
The goal of this paper is to deal with problem 2), if $\C{C}$ is the 
$C^*$-tensor category $\C{U}_G$ for a finite group $G$. 
The main result of this paper is that 
every subfactor $N\subset R$ of this kind is isomorphic 
to a subfactor $$R^G\subset (R\otimes \B{L}(\comp ^r))^H.$$
Here $R^G$ denotes the fixed point algebra under an outer action $\alpha $
of $G$ on $R$, $H$ is a subgroup of $G$, $\psi :H\pfeil \B{U}(\comp ^r)$ is a
projective representation of $H$ (satisfying a certain condition), and
$(R\otimes \B{L}(\comp ^r))^H$ denotes the fixed point algebra under an action
of $H$ where $H$ acts on $R$ by $\alpha \rest H$ and on $\B{L}(\comp ^r)$
via the conjugation with $\psi $. In the proof we have to generalize the method
used by M. Nakamura and Z. Takeda in order to determine the von Neumann
algebras between $R^G$ and $R$. The proof also uses induced
representations and the imprimitivity theorem. Unfortunately, it is possible
that different choices of $H$ and $\psi $ yield isomorphic subfactors.

This article is the abridged version of a part of the author's 
Habilitations\-schrift
\cite{Sch2}.

\section{Preliminaries} 
\subsection{$C^*$-tensor categories}
We suppose that all sets appearing in this article are small, i.e.
belong to a fixed universe.
Let $\C{C}$ be a $C^*$-tensor category. We assume that the objects of $\C{C}$
form a set. For two objects $\rho $ and 
$\sigma $ of $\C{C}$, the (tensor) product of $\rho $ and $\sigma $ is
denoted by $\rho \sigma $ and the space of morphisms with source
$\rho $ and target $\sigma $ by $(\rho ,\, \sigma )$. The product
of two morphisms $T\in (\rho ,\, \rho ')$ and
$S\in (\sigma ,\sigma ')$ is denoted by
$T\times S\in (\rho \sigma ,\, \rho '\sigma ')$. $\B{1}_{\sigma }$ denotes
the identity morphism of an object $\sigma $ of $\C{C}$.

A compact $C^*$-tensor
category is a strict $C^*$-tensor category for which every object $\sigma $
has a conjugate $\osigma $ and 
for which the space $(\iota ,\, \iota )$ is one dimensional
($\iota $ unit object). \cite{Sch4} contains a detailed definition.
A compact $C^*$-tensor category $\C{C}$
is called finite
if the number of the equivalence classes of the irreducible objects of
$\C{C}$ is finite. The statistical dimension of an object $\sigma $ of a
compact $C^*$-tensor category is denoted by $d(\sigma )$.

Let $\rho $ be an object of a finite $C^*$-tensor category
$\C{C}$. $\C{C}_{\rho }$ denotes the following
full $C^*$-tensor subcategory of $\C{C}$:
\newl The objects of $\C{C}_{\rho }$ form the smallest subset $\C{O}_{\rho }$
of the set of objects of $\C{C}$
with the following properties:
\begin{itemize}
\item[(a)] $\iota ,\,\rho \in \C{O}_{\rho }$.
\item[(b)] If $\tau \in \C{O}_{\rho }$, then every object equivalent to
$\tau $ and every subobject of $\tau $ belongs to $\C{O}_{\rho }$.
\item[(c)] If $\tau ,\,\phi \in \C{O}_{\rho }$, then $\ov{\tau }\in 
\C{O}_{\rho }$ and $\tau \phi \in \C{O}_{\rho }$.
\item[(d)] Finite direct sums of objects of $\C{O}_{\rho }$ are objects of
$\C{O}_{\rho }$.
\end{itemize}

Let $N\subset M$ be a an inclusion of $\zws $-factors
with finite Jones index \linebreak $[M:N]< \infty $ and with $N \not= M$  and 
let 
$$ N=M_{-1} \subset M=M_0\subset M_1 \subset M_2 \subset M_3 \subset \ldots $$
be the Jones tower for $N\subset M$. An $(N,N)$-bimodule is a Hilbert space
$\C{H}$ endowed with a normal left action $\lambda $ and a normal right action
$\rho $ of $N$ on $\C{H}$ such that $\lambda (N)$ and $\rho (N)$ commute.
 The $(N,N)$-bimodules which are
equivalent to a subbimodule of $_NL^2(M_k)_N$ for some $k\in \nat $ form
a compact $C^*$-tensor category $\C{B}_{N\subset M}$, 
where the product of objects is the $N$-tensor product $\tn $ (see 
\cite{Sch4}).
$\C{B}_{N\subset M}$ is the smallest $C^*$-tensor category of
$(N,N)$-bimodules containing the bimodule $_NL^2(M)_N$ as an object.

\subsection {$\zws $-subfactors defined by finite $C^*$-tensor categories
\label{S11}}
It is possible to construct a finite depth subfactor of the hyperfinite
$\zws $-factor for a given object $\sigma $ of a finite
 $C^*$-tensor category $\C{C}$ (see \cite{Sch4}, Section 3.3
compare also \cite{LoRo}. The construction seems to be known to
several people). Let $\osigma $ be an object conjugate to $\sigma $.
If $d(\sigma )>1$ the tower
\begin{equation} \hspace*{-8mm} \begin{array}{ccccccccc}
(\sigma ,\sigma ) & \subset &( \osigma \, \sigma   ,\, \osigma\, \sigma )&
\subset
 & (\sigma \, \osigma \, \sigma  ,\,\sigma \, \osigma  \, \sigma )   
&\subset &
(\osigma \,\sigma \, \osigma \, \sigma  ,\, \osigma \, \sigma \, 
\osigma \, \sigma  )
& \subset &\ldots \\
\cup & &\cup & &\cup & &\cup & & \\
(\iota , \iota )& \subset &(\osigma ,\,\osigma )& \subset &
(\sigma \, \osigma \, ,\, \sigma \, \osigma \, )  &\subset &  
(\sigma \, \osigma \, \sigma ,\, \sigma \, \osigma \, \sigma )
&\subset &\ldots 
\end{array} \label{T1} \end{equation}
of finite dimensional $*$-algebras fulfils the periodicity assumptions
used in H. Wenzl's subfactor construction (Theorem 1.5 in \cite{We}).
(The inclusion of two successive algebras in the same line
is given by $T\Pfeil \B{1}_{\sigma }\times T$ or
$T\Pfeil \B{1}_{\osigma }\times T$ and for two algebras in the same column,
it is given by
$T\Pfeil T\times \B{1}_{\sigma }$.)
Hence the tower defines a $\zws $-subfactor $A\subset B$, where
the union of the $*$-algebras in the upper line (resp. lower line) is
ultra-strongly dense in $B$ (resp. $A$). Obviously,
$B$ is the hyperfinite $\zws $-factor. We get $[B:A] = d(\sigma )^2$.
 Let
$$ A=B_{-1} \subset B=B_0\subset B_1 \subset B_2 \subset B_3 \subset \ldots $$
be the Jones tower for $A\subset B$. 
The standard invariant 
\[ \begin{array}{rcccccccc}
\comp \B{1}=A'\cap A &\subset &A' \cap B & \subset &A'\cap B_1&\subset 
& A'\cap B_2 &\subset & 
\ldots \\
 & &\cup & &\cup & &\cup & & \\
\comp \B{1}& =&B'\cap B & \subset &B' \cap B_1& \subset &B' \cap B_2 &\subset
&\ldots 
\end{array} \]  
of $A\subset B$ is equal to
\begin{equation} \hspace*{0mm}
\begin{array}{rcccccccc}
(\iota ,\iota ) &\subset &(\sigma ,\sigma ) & \subset &(\sigma \osigma , 
\sigma \osigma )
&\subset & (\sigma \osigma \sigma , 
\sigma \osigma \sigma ) &\subset &  
\ldots \\
 & &\cup & &\cup & &\cup  \\
& &(\iota , \iota ) & \subset & (\osigma ,\osigma )& \subset &(\osigma \sigma ,
\osigma \sigma ) &\subset & 
\ldots .
\end{array}   \label{T2}\end{equation}
Consequently, the subfactor $A\subset B$ has finite depth. The 
$C^*$-tensor category $\C{B}_{A\subset B}$ is equivalent to the full 
$C^*$-tensor subcategory $\C{C}_{\sigma \osigma }$ of $\C{C}$ (see 
\cite{Sch4}, Theorem 4.1).

\subsection{Subfactors of the hyperfinite $\zws $-factor 
belonging to a given finite 
$C^*$-tensor category $\C{C}$ \label{S12}}
Let a finite $C^*$-tensor category $\C{C}$ be given. We assume that 
$N\subsetneqq M=R$ is a subfactor of the hyperfinite $\zws $-factor $R$ with
finite Jones index such that $\C{B}_{N\subset M}$ is equivalent to
$\C{C}$. Let $A\subset B$ be the subfactor from Section ~\ref{S11}, where
the $C^*$-tensor category $\C{B}_{N\subset M}$ and the object $_NL^2(M)_N$
are used. Example 3.5 in \cite{Sch4} shows that $N\subset M_1$ and 
$A\subset B$ have the same standard invariant. So they are isomorphic according
to S. Popa's result \cite{P1} mentioned in the Introduction. If we replace
$\C{B}_{N\subset M}$ by the equivalent $C^*$-tensor category $\C{C}$ and
$_NL^2(M)_N$ by the corresponding object $\sigma $ of $\C{C}$, we obtain the
same subfactor $A\subset B$. We observe that $\sigma =\osigma $ and that the 
unit object $\iota $ is a subobject of $\sigma $. Furthermore we notice
$\C{B}_{N\subset M} =\C{B}_{N\subset M_1}$, which follows from
$L^2(M) \subset L^2(M_1)$ and $L^2(M_1)\cong L^2(M)\tn L^2(M)$.

If we want to find all subfactors $N\subsetneqq R$ of the 
hyperfinite $\zws $-factor $R$ for which the $C^*$-tensor category 
$\C{B}_{N\subset M}$ is equivalent to
$\C{C}$, we could proceed as follows:
\newl Regard the objects $\sigma $ of $\C{C}$ with the following properties:
\begin{itemize}
\item[(a)] $\sigma =\osigma $,
\item[(b)] $\iota $ is a proper subobject of $\sigma $,
\item[(c)] the $C^*$-tensor category $\C{C}_{\sigma }$ ($= 
\C{C}_{\sigma \osigma }$ by Property (a) and (b)) is equal to $\C{C}$.
\end{itemize}
Form the
subfactors $A\subset B$ for these objects $\sigma $ according to Section 
~\ref{S11} and determine the factors $C$ between $A$ and $B$ such that
the subfactor $A\subset B$ is isomorphic to the subfactor
$A\subset B(A,C)$ where $B(A,C)$ is the basic construction for 
the subfactor $A\subset C$. For every
subfactor $A\subset C$ of that kind the $C^*$-tensor category is in fact
equivalent to $\C{C}$. 

\subsection{Notation for groups and their representations \label{S13}}
The group algebra for a finite group $G$ is denoted by $\comp [G]$.
By a (unitary) projective representation of a finite group $G$
we mean a map $\pi $ from $G$ into the set $\B{U}(\C{H})$ of the unitary
operators of a Hilbert space $\C{H}$, for which there is a function
$c:G\times G \pfeil T:=\{z\in \comp :\, |z|=1\}$ such that 
\begin{equation} c(g,h)\, \pi (g)\, \pi (h) =  \pi (gh)\label{EE1} 
\end{equation} 
 holds for $g,h\in G$. 
The function $c:G\times G \pfeil T$
 is called a cocycle of $G$.
A projective representation $\pi $ of $G$ satisfying (~\ref{EE1}) is called
a $c$-projective representation of $G$.
The projective kernel 
$\proj \ker \pi $ of $\pi $ is $\proj \ker \pi :=
\{g\in G:\, \pi (g)\in \comp \B{1}\}$.

$\B{L}(\C{H})$ denotes the algebra of all continuous linear operators
on $\C{H}$.
We denote the action $x\in \B{L}(\C{H}) \Pfeil \pi (g) \, x \, \pi (g^{-1})$
of $G$ on $\B{L}(\C{H})$ by $\Ad \pi $.

Two projective representations $\pi _1:G\pfeil \B{U}(\C{H}_1)$ and
$\pi _2:G\pfeil \B{U}(\C{H}_2)$ are called equivalent, if there are a unitary
operator $U:\C{H}_1 \pfeil \C{H}_2$ and a function $\mu :G\pfeil T$
such that $U\, \pi _1(g)\, U^* = \mu (g)\, \pi _2(g)$ for every $g\in G$.

Let $H$ be a subgroup of $G$. We put $N(H):= \bigcap _{g\in G}gHg^{-1}$.
We use the notation $G/H$ for the space of the left cosets $kH$ in $G$.
We assume that a set $\cal V$ of representatives for the left cosets  
$kH,\, k\in G,$ is fixed
such that the neutral element $e$ belongs to $\cal V$.
For every $g\in G,\,\,
g=k(g) h(g)$ denotes the unique decomposition of $g$ such that 
$k(g)\in {\cal V}$ and $h(g)\in H$.
    
If $\pi :H\pfeil \B{U}(\C{K})$ is a finite dimensional representation
of $H$ we define the induced representation 
$\ind \pi :G\pfeil \B{U}(\C{K} \otimes \ell ^2(G/H))$ by
$$(\ind \pi )(g) \, (\xi \otimes \delta _{kH}) \, =
 \pi (h(gk))\, \xi \otimes \delta _{gkH} $$
for $\xi \in \C{K}, \, k\in \C{V}$ and $g\in G$.

\subsection{The subfactor $A\subset B$ for a finite group $G$ \label{S14}}
We consider the $C^*$-tensor category $\C{U}_G$ of the unitary finite
dimensional representations of a finite group $G$. The product of objects is 
the usual tensor product of representations, the object $\orho $ conjugate
to an object \linebreak
$\rho :G\pfeil \B{U}(\C{H})$ is the contragredient representation,
and $(\rho ,\rho )$ is the fixed point algebra $\B{L}(\C{H})$ under the
action $\Ad \rho $.

Recall that there exists an outer action of $G$ on the hyperfinite 
$\zws $-factor
$R$ and that two outer actions of $G$ on $R$ are conjugate.
The fixed point algebra $R^G$ under an outer action $\alpha $ of
$G$ is a $\zws $-factor.
If $\sigma :G\pfeil \B{U}(\C{H})$ is a finite dimensional representation
of $G$, $\alpha \otimes \Ad \sigma $ is an action of $G$ on
$R\otimes \B{L}(\C{H})$ and the fixed point algebra $(R\otimes \B{L}(\C{H}))^G$
is a $\zws $-factor. The subfactor $R^G\subset (R\otimes \B{L}(\C{H}))^G$
has index $(\dim \C{H})^2$.
If $\dim \C{H} >1$ the standard invariant of this subfactor is
equal to the standard
invariant of the subfactor $A\subset B$ from Section ~\ref{S11}
where the $C^*$-tensor category $\C{U}_G$ and the object $\sigma $ is used
(see \cite{W} or \cite{Sch}). According to S. Popa's result \cite{P1}
these subfactors are isomorphic.

\section{The Theorem}
\subsection{The subfactors $R^G\subset (R\otimes \B{L}(\comp ^r))^H$
\label{S15}}
As before we assume an outer action $\alpha $ of $G$ on $R$. Moreover, let
$\psi :H\pfeil \B{U}(\comp ^r)$ be a projective representation of
a subgroup $H$ of $G$. $\alpha \rest H \otimes \Ad \psi $ is an action of
$H$ on $R\otimes \B{L}(\comp ^r)$, let $(R\otimes \B{L}(\comp ^r))^H$ be the
fixed point algebra under this action. Then $R^G\subset
(R\otimes \B{L}(\comp ^r))^H$ is a $\zws $-subfactor with index
$[G:H]\cdot r^2$. In \cite{Sch} these subfactors are investigated and their
principal and dual principal graphs are computed. In \cite{Sch} only  
ordinary representations $\psi $ are used, but the transfer to projective
representations is easy. 

We review some facts from \cite{Sch}. The tensor product $\opsi \otimes \psi $
of the conjugate representation $\opsi $ and of $\psi $ is an ordinary 
representation of $H$ and the induced representation
$\ind (\opsi \otimes \psi )$ of $G$ is defined on the Hilbert space
$\C{V}:=\ov{\comp ^r}\otimes \comp ^r \otimes \ell ^2(G/H)$, where
$\ov{\comp ^r}$ is the Hilbert space dual to $\comp ^r$.
The $*$-algebra $\ell ^{\infty }(G/H)$ acts on $\ell ^2(G/H)$ 
by multiplication, the $*$-algebra 
$$C:=\comp \, id_{\ov{\comp ^r}}
\otimes \B{L}(\comp ^r)
 \otimes  \ell ^{\infty }(G/H) \cong \B{L}(\comp ^r)
 \otimes  \ell ^{\infty }(G/H) $$
of linear operators on $\C{V}$       
is invariant under the action $\Ad \ind (\opsi \otimes \psi )$ of $G$.
Let $e_{kH}\in \ell ^{\infty }(G/H)$ be the linear
operator on $\ell ^2(G/H)$ determined by $e_{kH}\, \delta _{lH} = 
\delta _{kH,lH}\, \delta _{kH}$ for $l\in \C{V}$. Every element $x$ of
$(R\otimes C)^G$ can be written uniquely as $x=\sum _{k\in \C{V}}
x_k \otimes e_{kH}$ where $x_k\in  R\otimes \B{L}(\comp ^r)$.
$x\in (R\otimes C)^G\Pfeil x_e$ defines an isomorphism $\iota $ from
$(R\otimes C)^G$ onto $(R\otimes \B{L}(\comp ^r))^H$, hence the subfactors
$R^G\subset (R\otimes C)^G$ and $R^G\subset
(R\otimes \B{L}(\comp ^r))^H$ are isomorphic. By using this result and
an invariance principle due to A. Wassermann, one finds that the subfactor
obtained from $R^G\subset
(R\otimes \B{L}(\comp ^r))^H$ by an application of the basic construction
is isomorphic to $R^G\subset
(R\otimes \B{L}(\C{V}))^G$ where $\Ad \ind (\opsi \otimes \psi )$ acts on
$\B{L}(\C{V})$. So we obtain 
$\C{B}_{R^G\subset
(R\otimes \B{L}(\comp ^r))^H} = \C{B}_{R^G\subset
(R\otimes \B{L}(\C{V}))^G}$. According to the Sections ~\ref{S11}, ~\ref{S12}
and ~\ref{S14}, the $C^*$-tensor category 
$\C{B}_{R^G\subset
(R\otimes \B{L}(\comp ^r))^H}$ is equivalent to $(\C{U}_G)_
{\ind (\opsi \otimes \psi ) }$.
The kernel $K$ of $\ind (\opsi \otimes \psi )$ is equal to the 
projective kernel $\proj \ker \psi \rest N(H)$ of the restriction
of $\psi $ onto $N(H)$. 

By applying Theorem 27.39 in \cite{HR} and observing 
$\ov{\ind (\opsi \otimes \psi )}=\linebreak \ind (\opsi \otimes \psi )$,
we obtain that every irreducible representation of $G$ is
contained in the $n$-fold tensor product 
$\bigl(\ind (\opsi \otimes \psi )\bigr)^{\otimes n}$ for 
some $n\in \nat $ if and only if
$K=\{e\}$. So $\C{B}_{R^G\subset
(R\otimes \B{L}(\comp ^r))^H}$ is equivalent to $\C{U}_G$ if and only if
$\proj \ker \psi \rest N(H) = \{e\}$.

The subfactor $R^G\subset
(R\otimes \B{L}(\comp ^r))^H$ is irreducible if and only if
$\psi $ is irreducible. We even have the more general result
$$(R^G)'\cap (R\otimes \B{L}(\comp ^r))^H  =\comp \B{1} \otimes \psi (H)'.$$
The inclusion $'\subset '$ is obvious, the other inclusion follows
from 
$$(R^G)'\, \cap \, \Bigr(R\otimes \B{L}(\comp ^r)\Bigr)^H \subset
\Bigl(\bigl((R^G)'\cap R\bigr) \otimes \B{L}(\comp ^r)\Bigr)^H$$

\abs Now let us formulate the main result of this article:

\begin{thm}
Let $G$ be a finite group and
let $N\subset M=R$ be a 
subfactor of the hyperfinite $\zws $-factor $R$ with finite Jones index
such that the $C^*$-tensor category 
$\C{B}_{N\subset M}$ is equivalent to $\C{U}_G$.
Then there are a subgroup $H$ of $G$, a projective finite dimensional
representation $\psi :H\pfeil \B{U}(\comp ^r)$ ($r\in \nat $) of $H$ satisfying
$\hbox{\rm proj }\ker \psi \rest N(H)\, =\{e\}$
and an outer action $\alpha $ of $G$ on the hyperfinite $\zws $-factor $R$
such that $N\subset M$ is isomorphic to 
\begin{equation} R^G \subset (R\otimes \B{L}(\comp ^r))^H.  \label{E1} 
\end{equation}
\label{Th} \end{thm} 

\subsection{Remarks: \label{S17}}
(1) Let $R^G \subset (R\otimes \B{L}(\comp ^r))^H$ and 
$R^{\tilde{G}} \subset (R\otimes \B{L}(\comp ^{\tilde{r}}))^{\tilde{H}}$
be two subfactors as in Theorem ~\ref{Th} and let $\tilde{\psi }:\tilde{H}
\pfeil \B{U}(\comp ^{\tilde{r}})$ denote the projective
representation used in the second subfactor.

If there is an isomorphism $\gamma $ from $G$ onto $\tilde{G}$ such that
$\gamma (H)=\tilde{H}$ and such that the projective representation
$\tilde{\psi }\circ \gamma $ of $H$ is equivalent to the representation
$\psi $, then the subfactors are isomorphic.

The converse is not true, as the following easy example shows:
\newl Let $G=\ganz _2 \times \ganz _2$, $G=H$ and $\psi:
\ganz _2 \times \ganz _2\pfeil \B{U}(\comp ^2)$ be the projective
representation of $\ganz _2 \times \ganz _2$ given by the Pauli matrices:
\begin{eqnarray*}
\psi (0,0) =\left( \begin{array}{cc} 1& 0 \\ 0&1 \end{array} \right),
& \psi (1,0) =\left( \begin{array}{cc} 0& 1 \\ 1&0 \end{array} \right), \\
\psi (0,1) =\left( \begin{array}{cc} 0& -i \\ i&0 \end{array} \right),
& \psi (1,1) =\left( \begin{array}{cc} 1& 0 \\ 0&-1 \end{array} \right).
\end{eqnarray*}
The subfactor $R^G \subset (R\otimes \B{L}(\comp ^2))^G$ given by these
data is an irreducible subfactor with index 4 and depth 2. (The principle
graph can be computed as in \cite{Sch}. Observe that every irreducible 
projective representation with the same cocycle is equivalent to $\psi $.)
All irreducible subfactors of the hyperfinite $\zws $-factor with index
$4$ and depth $2$ are isomorphic to $R^G\subset R$ or $R^{\ganz _4}\subset R$.
The $C^*$-tensor category $\C{B}_{R^G \subset (R\otimes \B{L}(\comp ^2))^G}$
is equivalent to $\C{U}_G$ and $\C{B}_{R^{\ganz _4}\subset R}$ to
$\C{U}_{\ganz _4}$. Moreover, $\C{U}_G$ and $\C{U}_{\ganz _4}$ are not 
equivalent. So the subfactor $R^G \subset (R\otimes \B{L}(\comp ^2))^G$
is isomorphic to $R^G\subset R$.

\abs
\noindent (2) The $C^*$-tensor categories associated with the subfactors
dual to the subfactors (~\ref{E1}) are investigated in \cite{Sch2}. If 
$\psi $ is an ordinary representation they can be described by the
$(H,H)$-vector bundles over $G$ introduced in \cite{KY} by H. Kosaki and
S. Yamagami. The general case requires a slight generalization.

\abs
\noindent (3) It would be interesting to find an answer for the following
question:
 
\abs 
{\em Does there exist non-isomorphic finite groups $G_1$ and $G_2$
for which the $C^*$-tensor categories $\C{U}_{G_1}$ and $\C{U}_{G_2}$  
are equivalent?}

\abs 
\noindent Observe that neither the Tannaka-Krein theorem nor
S.~Doplicher's and J.E.~Roberts' results concerning compact groups
and $C^*$-tensor categories with a symmetry (see \cite{DR}) imply the
answer no. There is an example of a finite group $G_0$ and a finite dimensional
non-commutative Hopf-$*$-algebra $A$ such that $\C{U}_{G_0}$ and the
$C^*$-tensor category $\C{U}_A$ of all finite dimensional unitary
corepresentations of $A$ are equivalent (see \cite{Sch2}). We intend
to deal with this example in a forthcoming paper.

\section{The Proof of the Theorem}
Let $G\not= \{e\}$ and let $N\subset M=R$ be a subfactor of the 
hyperfinite $\zws $-factor
$R$ with finite index such that the
$C^*$-tensor category $\C{B}_{N\subset M}$ is equivalent to $\C{U}_G$.
Section ~\ref{S12} and Section ~\ref{S14} show that
there is a finite dimensional representation $\sigma :G\pfeil \B{U}(\C{H})$
of $G$ such that the subfactor $N\subset M_1$ is isomorphic
to the subfactor 
$$R^G \subset (R\otimes \B{L}(\C{H}))^G$$
where $G$ acts on $R$ by an outer action $\alpha $ and on $R\otimes
\B{L}(\C{H})$ by $\alpha \otimes \Ad \sigma $.

As $M$ is a factor between $N$ and $M_1$, there is a factor $P$ between
$R^G$ and $ (R\otimes \B{L}(\C{H}))^G$ such that $N\subset M$ is isomorphic to
$R^G\subset P$. We intend to determine the factors $P$ 
between $R^G$ and $(R\otimes \B{L}(\C{H}))^G$.

In \cite{NaTa1} and \cite{NaTa2} M. Nakamura and
Z. Takeda showed that every von Neumann algebra between $R^G$ and $R$ 
is equal to $R^H$ for some subgroup $H$ of $G$. 
Some parts of our subsequent considerations may be regarded as a 
generalization
of the methods used there.
We will investigate the von Neumann algebras between 
$((R\otimes \B{L}(\C{H}))^G)'$ and $(R^G)'$ where the commutants are formed 
in $L^2(R) \otimes \C{H}$. $R \otimes \B{L}(\C{H})$ acts on 
$L^2(R)\otimes \C{H}$ by $(r\otimes x)\, \ov{s} \otimes \xi =
\ov{rs} \otimes x\, \xi $ for $r\in R,\, x\in \B{L}(\C{H}),\, s\in R$ and
$\xi \in \C{H}$. ($\ov{s} $ denotes the element $s\in R$, if it is regarded
as an element of $L^2(R)$.)

The commutant $R'$ of $R$ in $L^2(R)$ is $S=\{\rho (r):\, r\in R\}$ 
($\rho $ right multipli\-cation). 
$\beta (g)\, \rho (r)\, = 
\rho (\alpha (g)r)$ for $g\in G$ and $r\in R$ defines an outer action
$\beta $ of $G$
on $S$.

The commutant $(R^G)'$ in 
$L^2(R)$ is generated by $S$ and $\{u_g:\, g\in G\}$ 
where $u_g$ is the unitary operator in $L^2(R)$ defined by
$u_g \,\ov{r} =\ov{\alpha (g)r}$ for $r\in R$, so the commutant can be 
identified with the crossed product $S \rtimes _{\beta } G$, 
which is a $\zws $-factor.  
It follows that the commutant $(R^G)'$ in $L^2(R) \otimes \C{H}$ is equal to
$(S\kreuz G) \otimes \B{L}(\C{H})$.

We determine the commutant $((R\otimes \B{L}(\C{H}))^G)'$ in
$L^2(R) \otimes \C{H}$: 
Obviously there is
a homomorphism $j$ from $S\kreuz G$ into $((R\otimes \B{L}(\C{H}))^G)'$
given by \linebreak
$j(x) =x \otimes \B{1}$
for $x\in S$ and $j(u_g) = u_g \otimes \sigma (g)$ for $g\in G$.
$j$ is injective and normal and the 
image $T:=j(S\kreuz G)$ is a $\zws $-factor contained in \newl
$\bigl((R\otimes \B{L}(\C{H}))^G\bigr)'$.
In fact we have 
$$T= ((R\otimes \B{L}(\C{H}))^G)'.$$
In order to show the last equation it suffices
 to verify \linebreak
$[T:j(S)]\, =[((R\otimes \B{L}(\C{H})^G))':j(S)]$.
We have  $[T:j(S)] = [S\kreuz G:S] = |G|$ as well as
\begin{eqnarray*}
 \bigl[ \bigl((R\otimes \B{L}(\C{H}))^G\bigr)':j(S)\bigr] &=&
\bigl[ \bigl((R\otimes \B{L}(\C{H}))^G \bigr)':(R\otimes \B{L}(\C{H}))'
\bigr]\, = \\
\hbox{$[ R\otimes \B{L}(\C{H}):(R\otimes \B{L}(\C{H}))^G]$} &=& |G|. 
\end{eqnarray*}
It is possible to derive the last equation from
$$\displaylines{ (\dim \C{H})^2\cdot |G| \,= [R\otimes \B{L}(\C{H}):R] \cdot 
[R:R^G]\, =
[R\otimes \B{L}(\C{H}):R^G]= \hfill \cr 
\hfill = [R\otimes \B{L}(\C{H}):(R\otimes \B{L}(\C{H}))^G] \cdot
[(R\otimes \B{L}(\C{H}))^G:R^G]    }$$ 
and
\begin{equation} [(R\otimes \B{L}(\C{H}))^G:R^G]  = 
(\dim \C{H})^2\label{E31Inv} \end{equation}
(see \cite{Sch} or \cite{W} as to a proof of equation (~\ref{E31Inv})).

\begin{lem} (i) 
Let $Q$ be the finite dimensional $*$-subalgebra of $S\kreuz G \otimes
\B{L}(\C{H})$ generated by $\{u_g \otimes x:\,g\in G,\, x\in \B{L}(\C{H})\}$,
which obviously is isomorphic to $\comp [G] \otimes \B{L}(\C{H})$. 
Let $\C{D}$ denote the set of all $*$-subalgebras $D$ of $Q$ with the same unit
like $Q$ and with the following property:
There exist $\comp $-vector spaces $B_g\subset \B{L}(\C{H}),\, g\in G,$
such that $D$ is the direct sum
$\bigoplus _{g\in G} u_g \otimes B_g$ as a $\comp $-vector space 
($B_g=0$ is permitted).
\newl The map $F\Pfeil F\cap Q$ is a bijection between the set
$$\C{F}:=\{F \hbox{a von Neumann algebra: }   
S \subset F \subset S\kreuz G \otimes \B{L}(\C{H})\}$$ and $\C{D}$.
\newl (ii)
If $F$ is a von Neumann algebra belonging to $\C{F}$ and \linebreak
$E_F:S\kreuz G \otimes \B{L}(\C{H}) \pfeil F$ denotes the conditional 
expectation onto $F$ corresponding to the unique normalized trace $\tr $ of
$S\kreuz G \otimes \B{L}(\C{H})$, then
for any $s\in S,\, g\in G$ and $x\in
\B{L}(\C{H})$ there exists an operator $a\in \B{L}(\C{H})$
such that
$E_F(su_g \otimes x) =su_g \otimes a$.
\label{L51} \end{lem}

\B{Proof:\ }Let $F$ be a von Neumann algebra in $\C{F}$, 
let $(v_1,\ldots ,v_m)$
be an orthonormal base of $\C{H}$ and let
$\epsilon _{ij} \in \B{L}(\C{H})$ ($i,j=
1,\ldots ,m$) be defined by $\epsilon _{ij} v_k =\delta _{j,k}\, v_i$ for
$k=1,\ldots ,m$.
For $x\in 
\B{L}(\C{H})$ and $g\in G$ there are unique elements $s_{h,ij}\in S$
($h\in G,$   $i,j=1,\ldots ,m$) such that 
\begin{equation} E_F(u_g \otimes x)\, =\sum _{h\in G} 
\sum _{i,j=1}^m s_{h,ij}\, 
     u_h \otimes \epsilon _{ij}. \label{E5E} \end{equation}
$$\displaylines{ E_F(u_g \otimes x)\cdot (n\otimes \B{1})\, =
   E_F(u_g n\otimes x)\, =E_F(u_gnu_g^{-1}u_g \otimes x)\, = \hfill \cr
  \hfill  E_F(\alpha (g)n\cdot u_g \otimes x) \,= 
  (\alpha (g)n\otimes \B{1})\cdot E_F(u_g \otimes x)      }$$
holds for every $n\in S$, hence
$$\sum _{h,i,j} s_{h,ij}\, u_h\, n \otimes \epsilon _{ij} \,= \,
   (\alpha (g)n\otimes \B{1})\, \sum _{h,i,j} s_{h,ij}\, u_h \otimes 
\epsilon _{ij}. $$
Using $u_h\cdot n\, =\alpha (h)n\cdot u_h$ we obtain 
$s_{h,ij}\cdot \alpha (h)n\, = \, \alpha (g)n\cdot s_{h,ij}$ for 
$h\in G,\, n\in S$ and $i,j=1,\ldots ,m$. Replacing $n$ by $\alpha (g^{-1})
n$ we get
$$s_{h,ij}\cdot \alpha (hg^{-1})n\, = \, n\cdot s_{h,ij}$$
for every $n\in S$.
If $h=g$ this relation means that $s_{g,ij}$ belongs to 
\linebreak $S\cap S'=\comp \B{1}$. Since 
the action $\alpha $ is outer, for
$h\in G\setminus \{e\}$ there is 
no element $s\in S\setminus \{0\}$ such that
$s\cdot \alpha (h)n= n\cdot s$ holds for every $n\in S$.
This implies $s_{h,ij}=0$ for $h\neq g$. By summarizing the considerations
from above we obtain $E_F(u_g \otimes x)=u_g\otimes a$ for a suitable
operator $a\in \B{L}(\C{H})$. Now (ii) follows immediately.

Moreover, we get $E_F(Q)\subset Q$ and $E_F(Q)=F\cap Q$.
If $\sum _{g\in G}u_g \otimes x_g\in F\cap Q$, then
$$\sum _{g\in G}u_g \otimes x_g\, = 
E_F\Bigl(\sum _{g\in G}u_g \otimes x_g\Bigr)\, =
  \sum _{g\in G}E_F(u_g \otimes x_g). $$
$E_F(u_g\otimes x_g)\in u_g\otimes \B{L}(\C{H})$ implies
$u_g\otimes x_g =E_F(u_g\otimes x_g)\in F\cap Q$ for every $g\in G$. 
Hence $F\cap Q$ is a direct sum $\bigoplus _{g\in G}u_g\otimes B_g$.

The assignment $E_F\Pfeil E_F \rest Q$ is injective,
furthermore $E_F\rest Q$ is the
conditional expectation from $Q$ onto $F\cap Q$ with respect to the trace 
\linebreak $\tr \rest Q$ such that $E_F\rest Q$ is determined by its image
$F\cap Q$. Hence the map \linebreak
$F\in \C{F}\Pfeil F\cap Q\in \C{D}$ is injective.

We have to show that this map is surjective. Let 
$C=\bigoplus _{g\in G} u_g\otimes B_g$ be a $*$-algebra in $\C{D}$.
For $g\in G$ let
$(b_g^1,\ldots ,b_g^{n(g)})$ ($n(g)\in \{0,1,\dots ,m^2\}$) be a base of the
$\comp $-vector space $B_g$. 
$$F:=\Bigl\{\sum _{g\in G} \sum_{j=1}^{n(g)} s_{g,j}\, u_g \otimes b_g^j:\,
 s_{g,j}\in S \hbox{ for $g\in G$ and $j=1,\ldots ,n(g)$} \Bigr\} $$
is a $*$-algebra satisfying $S\subset F \subset S\kreuz G \otimes
\B{L}(\C{H})$ and $F\cap Q =C$. It is not difficult to prove that
$F$ really is a von Neumann algebra.  One easily sees that the proof of this
fact is not necessary for the proof of Theorem ~\ref{Th}, so we omit it.
\blacksquare

\begin{lem} Via $F\Pfeil F\cap (\comp \B{1} \otimes  
 \B{L}(\C{H}) )$, we have a bijective correspondence between the von Neumann
algebras $F$ satisfying \linebreak
$T\subset F \subset S\kreuz G \otimes \B{L}(\C{H})$
and the $*$-subalgebras $B$ of $\comp \B{1} \otimes  
 \B{L}(\C{H}) \cong \B{L}(\C{H})$, which have the same unit as 
$\B{L}(\C{H})$ and are invariant under the action $\Ad \sigma $ of $G$. 
The inverse map is 
$$B\Pfeil F_B:=\Bigl\{\sum _{g\in G} \sum _{j=1}^k s_{g,j}u_g \otimes
  \sigma (g) \, b^j: \, s_{g,j}\in S \hbox{ for $g\in G,\, j=1,\ldots ,k$} 
   \Bigr\},$$
where $(b^1,\ldots ,b^k)$ is a base of $B$.
\label{L52} \end{lem}

\B{Proof:\ } Let $F$ be a von Neumann algebra between $T$ and
$S\kreuz G \otimes \B{L}(\C{H})$ and
$F\cap Q =$ 
 $\bigoplus _{g\in G} u_g \otimes B_g$ be the decomposition of $F\cap Q$ as
in Lemma ~\ref{L51}. We put
$B:=B_e$. Obviously $B$ is a $*$-algebra. We get
\begin{equation}
 \sigma (g)B= B_g \qquad \hbox{for every $g\in G$,}  \label{E5s} \end{equation}
as the following conclusion shows:
\begin{eqnarray*}
x\in B_g & \Longleftrightarrow & E_F(u_g\otimes x) =u_g\otimes x \qquad
  \Longleftrightarrow \qquad \qquad 
\hbox{(note $u_g \otimes \sigma (g) \in T$)}  \\
 &\Longleftrightarrow &
  (u_g \otimes \sigma (g)) \cdot E_F(\B{1}\otimes \sigma (g^{-1})\, x) \, =
(u_g \otimes \sigma (g)) \cdot (\B{1}\otimes \sigma (g^{-1})\, x) \\
 &\Longleftrightarrow &  E_F(\B{1}\otimes \sigma (g^{-1})\, x) \, =
 \B{1}\otimes \sigma (g^{-1})x
 \qquad \Longleftrightarrow \\
 &\Longleftrightarrow &
\sigma (g^{-1})x\in B.  \end{eqnarray*} 
The injectivity of $F\Pfeil F\cap Q$ implies that the assignment
$F\Pfeil B$ is injective, too.

If we multiply $u_g \otimes \sigma (g)$ from right instead of from left we get
\begin{equation}
 B \,\sigma (g)= B_g \qquad \hbox{for $g\in G$.}  \label{E5s2} \end{equation}
The relations (~\ref{E5s}) and (~\ref{E5s2}) imply 
$\sigma (g)\,  B\, \sigma (g^{-1})  \,=B$ for every $g\in G$.

Conversely, let $B$ be a $*$-subalgebra of $\B{L}(\C{H})$ with the same unit 
invariant under the action $\Ad \sigma $. It is easy to check that
$D=\bigoplus _{g\in G} u_g \otimes \sigma (g)B$ is  a $*$-subalgebra of 
$Q$. According to the proof of Lemma ~\ref{L51}, $F_B$ is the 
von Neumann algebra in $L^2(R)\otimes \C{H}$ associated with $D$. Clearly,
\linebreak $F_B\cap  (\comp \B{1} \otimes \B{L}(\C{H})) \, =B$ 
and $T\subset F_B\subset S\kreuz G \otimes \B{L}(\C{H})$ hold.
\blacksquare

\begin{lem} The center $Z_F$ of a von Neumann algebra $F$ between
$T$ and $S\kreuz G \otimes \B{L}(\C{H})$ consists of those elements 
of the center $Z_B$ of \newl $B:=F \cap (\comp \B{1} \otimes \B{L}(\C{H}))$
which are invariant under $\Ad \sigma $. \label{L53} \end{lem}
 
\B{Proof:\ }
$Z_F =F'\cap F \subset (S'\otimes \B{L}(\C{H})) \cap 
(S\kreuz G \otimes \B{L}(\C{H}))
=\comp \B{1}\otimes \B{L}(\C{H})$ 
(observe $S'\cap S\kreuz G =\comp \B{1}$) and consequently
$Z_F \subset B$. $B\subset F$ implies $Z_F \subset Z_B$.
Since the elements of $Z_F$ commute with $u_g\otimes \sigma (g)$ for $g\in G$,
we obtain $Z_F \subset \{x\in Z_B:\, \sigma (g)\, x\, \sigma (g^{-1}) \,=x\}$. 

Conversely, every element $x$ of $Z_B$ invariant under $\Ad \sigma $
commutes with $S\otimes \comp \B{1}$, $B$ and $u_g\otimes \sigma (g)$ for 
$g\in G$. Since the von Neumann algebra $F$ is generated by these elements,
$x$ belongs to $Z_F$. \blacksquare

\begin{lem} Let $F$ be a factor between $T$ and 
$ S\kreuz G \otimes \B{L}(\C{H})$ and let 
$B:=\linebreak F\cap (\comp \B{1} \otimes
\B{L}(\C{H}))$. There are a subgroup $H$ of $G$, projective
representations $\rho :H\pfeil \B{U}(\comp ^d)$ and
$\psi :H\pfeil \B{U}(\comp ^r)$ and a unitary operator
$U:\C{H} \pfeil \comp ^d \otimes \comp ^r \otimes \ell ^2(G/H)$
such that
$\rho \otimes \psi $ is an ordinary representation of $H$,
\begin{eqnarray}
U \, \sigma (g)\, U^* &=& \ind (\rho \otimes \psi )(g) \qquad
   \hbox{for $g\in G$ \ \ \ \ and} \label{E5U1} \\
U\, B\,  U^*  &=& \B{L}(\comp ^d) \otimes \comp \, id_{\comp ^r} \otimes
   \ell ^{\infty }(G/H).     \label{E5U2} \end{eqnarray}
(In particular, if
$c:H\times H \pfeil T$ is the cocycle for $\psi $, so \vspace{1pt} 
$\ov{c}(h_1,h_2) =\ov{c(h_1,h_2)}\, (h_1,\, h_2\in H )$ is the cocycle
for $\rho $.)
\label{L54} \end{lem}

\B{Proof:\ }The proof of Lemma ~\ref{L52} showed 
$\sigma (g)\, B\, \sigma (g^{-1})
=B$ for every $g\in G$. Since every automorphism of an algebra maps the center 
of this algebra onto itself, $\Ad \sigma $ leaves the center $Z_B$ of $B$
invariant. So the group $G$ acts on the set $\{p_1,\ldots ,p_l\}$ of
the minimal projections of the center $Z_B$ by
$g.p_j \,=\sigma (g) p_j \sigma (g)^{-1}$. This action is transitive, otherwise
there would exist a projection $p\neq 0,\B{1}$ of $Z_B$ 
satisfying $p = g. p$ for every $g\in G$. Hence
$p$ would belong to the center $Z_F$ of $F$ by Lemma ~\ref{L53},
which is a contradiction to the assumption that $F$ is a factor.

So we are able to apply the imprimitivity theorem. There are a subgroup $H$, 
an ordinary representation
$\pi :H\pfeil \B{U}(\C{K})$ of $H$ and a unitary operator 
$U:\C{H} \pfeil \C{K} \otimes \ell ^2(G/H)$ such that
\begin{eqnarray*}
U\, \sigma (g)\,U^* &=& \ind \pi (g) \qquad \hbox{for $g\in G$ and}\\
\{U\,p_j\, U^*:\, j=1,\ldots ,l\} &=&\{ id_{\C{K}} \otimes e_{kH}:\,
 k\in \C{V} \}    \end{eqnarray*} 
(see Section ~\ref{S15} for the definition of $e_{kH}$).
Without loss of generality we assume 
$U\, p_1 \,U^* =id_{\C{K}}\otimes e_H$. 

From now on we will identify $\C{H}$ and $\C{K}\otimes \ell ^2(G/H)$ 
as well as $\sigma $ and $\ind \pi $. The finite dimensional factor
$Bp_1$ may be identified with a subfactor of $\B{L}(\C{K}) \otimes \comp \, e_H
\cong \B{L}(\C{K})$. There are finite dimensional Hilbert spaces 
$\comp ^d$ and $\comp ^r$ and a unitary operator
$V:\C{K} \pfeil \comp ^d \otimes \comp ^r$ such that
$$V\, Bp_1\, V^* \, =\, \B{L}(\comp ^d)\otimes \comp \, id_{\comp ^r}
 $$ holds.

Using $V$ we identify $\C{K}$ and $\comp ^d \otimes \comp ^r$. We fix
an element $h\in H$.  
From $\sigma (h)\,p_1 \,\sigma (h)^* =p_1$ we conclude that
$$\displaylines{x = y\otimes \B{1} \otimes e_H \in Bp_1 = 
\B{L}(\comp ^d)\otimes \comp \B{1}
  \otimes \comp e_H
\Pfeil \sigma (h)\,x\,\sigma (h)^* = \hfill\cr \hfill
\pi (h) (y\otimes \B{1}) \pi (h)^* \otimes e_H   }$$
is an automorphism
of $Bp_1$. This automorphism is inner, 
hence there is an operator $\rho (h)\in \B{U}(\comp ^d)$ such that
\begin{equation} \pi (h) (y\otimes id_{\comp ^r}) \pi (h)^* \,=
\, \rho (h)\, y \rho (h)^* \otimes id_{\comp ^r} \label{E5r} \end{equation}
for $y\in \B{L}(\comp ^d)$. 
Let $h_1$ and $h_2$ be elements of $H$. For $y \in \B{L}(\comp ^d)$ we have
$$\displaylines{
\rho (h_1 h_2)\, y\,\rho (h_1 h_2)^* \otimes \B{1} \, =\,
\pi(h_1 h_2) \, (y\otimes \B{1}) \pi(h_1h_2)^* \,= \hfill\cr
\hfill \pi (h_1) 
\Bigl(\pi(h_2) \, (y\otimes \B{1})\, \pi(h_2)^*\Bigr) \pi (h_1)^*\,=
\rho (h_1) \rho (h_2)\, y\,\rho (h_2)^* \rho (h_1)^* \otimes \B{1}.   }$$
Hence $\rho (h_1h_2)^{-1} \rho (h_1)\rho (h_2)
\in \B{L}(\comp ^d)' =\comp \B{1}.$
So there is a $c(h_1,h_2)\in T$ such that \vspace{1pt}
$\ov{c(h_1,h_2)} \rho (h_1)\, \rho (h_2) =\rho (h_1h_2)$, hence
$\rho :H\pfeil \B{U}(\comp ^d), \,h\Pfeil \rho (h),$
is a projective representation. Relation (~\ref{E5r}) implies that
$(\rho (h) \otimes \B{1})^{-1}\pi (h)$ is a unitary operator in
$(\B{L}(\comp ^d)\otimes \comp \, id_{\comp ^r})'=
 \comp \, id_{\comp ^d} \otimes \B{L}(\comp ^r)$. Therefore there is a 
unique unitary operator $\psi (h)\in \B{U}(\comp ^r)$ such that 
$\pi (h)=\rho (h) \otimes \psi (h)$. Since $\pi $ is an
ordinary representation of $H$, 
$\psi :H\pfeil \B{U}(\comp ^r), \,h\Pfeil \psi (h),$
is a projective representation of $H$ with cocycle $c$. 

For every $j=1,\ldots ,l$ there is a $g\in G$ and a $k\in \C{V}$ such that
$p_j \, = g.p_1\, =id_{\C{K}} \otimes e_{kH}$. Hence
\begin{align}
 B\, p_j\, &= \sigma (g) \, Bp_1\, \sigma (g)^{-1}\, \notag\\
&= \ind (\rho \otimes \psi )(g)\, \Bigl(\B{L}(\comp ^d) \otimes \comp \, \B{1}
\otimes \comp \, e_H \Bigr)\, \ind (\rho \otimes \psi )(g)^{-1} \, \notag \\
&  =\, \B{L}(\comp ^d) \otimes \comp \, \B{1}
\otimes \comp \, e_{kH}, \notag \end{align}
and Equation (~\ref{E5U2}) follows. 
\blacksquare

\abs
Every factor between $R^G$ and $(R\otimes \B{L}(\C{H}))^G$ is the
commutant $F'$ of a factor $F$ between 
$T=\Bigl((R\otimes \B{L}(\C{H}))^G\Bigr)'$
and $(R^G)'=S\kreuz G \otimes \B{L}(\C{H})$.
Let $F$ be defined as in Lemma
~\ref{L54}. By applying $U$ we identify $\C{H}$ and $\comp ^d \otimes
\comp ^r \otimes \ell ^2(G/H)$ again.

Let $C$ be the $*$-algebra $\comp \, id_{\comp ^d} \otimes \B{L}(\comp ^r)
\otimes \ell ^{\infty }(G/H)$. We get
$$ (S\otimes B)' \, =\, S'\otimes B' \, =\, R \otimes
\Bigl(\B{L}(\comp ^d)\otimes \comp \B{1} \otimes 
\ell ^{\infty }(G/H)\Bigr)' \,=
\, R \otimes C     $$
and 
$$F' =\{x\in R\otimes C: \, \hbox{$x$ and $u_g \otimes \sigma (g)$ commute
  for every $g\in G$}  \} \, = (R\otimes C)^G. $$

So every subfactor $N\subset M=R$ with finite index for which the $C^*$-tensor
category $\C{B}_{N\subset M}$ is equivalent to $\C{U}_G$ is isomorphic
to a subfactor $R^G\subset (R\otimes C)^G$. The arguments from Section 
~\ref{S15} show that this subfactor is equivalent to 
$R^G\subset (R\otimes \B{L}(\comp ^r))^H$ and that
the restriction of the projective kernel of $\psi $ to $N(H)$ has to be
$\{e\}$.

\end{document}